\documentclass[reqno,11pt]{amsart}
\usepackage{amssymb}

\usepackage{amsmath}
\usepackage{times}
\usepackage{graphics}
\usepackage{hyperref}

\newtheorem{theorem}{Theorem}[section]
\newtheorem{lemma}[theorem]{Lemma}
\newtheorem{proposition}[theorem]{Proposition}

\newtheorem{corollary}[theorem]{Corollary}

\begin{document}
\title{On Fluctuations of Riemann's Zeta Zeros}
\author{V. Kargin }
\thanks{Statistical Laboratory, Department of Mathematics, University of
Cambridge, Cambridge, UK. e-mail: v.kargin@statslab.cam.ac.uk}
\date{September 2012}
\maketitle

\begin{center}
\textbf{Abstract}
\end{center}

\begin{quotation}
It is shown that the normalized fluctuations of Riemann's zeta zeros around
their predicted locations follow the Gaussian law. It is also shown that
fluctuations of two zeros, $\gamma _{k}$ and $\gamma _{k+x},$ with $x\sim
\left( \log k\right) ^{\beta },$ $\beta >0,$ for large $k$ follow the
two-variate Gaussian distribution with correlation $\left( 1-\beta \right)
_{+}.$
\end{quotation}

\section{Introduction}

This paper is concerned with the statistical properties of the Riemann zeta
function zeros. This subject originated in 1944, when Selberg \cite%
{selberg44} showed that the number of zeros in a sufficiently long interval
on the critical line can be described by the Gaussian law (see also \cite%
{selberg46}, \cite{ghosh83}, \cite{joyner86}, \cite{hughes_nikeghbali_yor08}%
). In the 1970s, Montgomery and Dyson discovered the remarkable fact that
the spacings between the zeta zeros resemble the spacings between the
eigenvalues of random Hermitian Gaussian matrices. This resemblance was
substantiated analytically by Montgomery \cite{montgomery73} and supported
numerically by Odlyzko \cite{odlyzko87} (see also \cite{rudnick_sarnak96}, %
\cite{coram_diaconis03}). The connection between zeta zeros and random
matrix eigenvalues drew much attention, as can be seen for example from
review papers in \cite{mezzadri_snaith05}. Recently, Bourgade \cite%
{bourgade10} supported this connection by showing that at the mesoscopic
level Riemann's zeros have correlations previously found by Diaconis and
Evans \cite{diaconis_evans01} for eigenvalues of unitary random matrices.
(See also \cite{bourgade_kuan12}.)

The motivation for our study comes from a paper by Gustavsson \cite%
{gustavsson05}, who showed that eigenvalues of random Hermitian matrices
fluctuate according to the Gaussian law. Our goal is to investigate the
statistical fluctuations of Riemann's zeros around their predicted positions
and to show that these fluctuations also follow the Gaussian law.

We denote the non-trivial zeros of Riemann's zeta function by $\beta
_{k}+i\gamma _{k}.$ (We do not assume Riemann's hypothesis in this paper.)
We consider only zeros with $\beta _{k}\geq 1/2$ and positive imaginary
part, $\gamma _{k}>0$, and order them so that the imaginary part is
non-decreasing, $\gamma _{1}\leq \gamma _{2}\leq \ldots .$

Let $\mathcal{N}(T)$ denote the number of zeros with the imaginary part
strictly between $0$ and $T$. If there is a zero with imaginary part equal
to $T$, then we count this zero as $1/2$.

Define 
\begin{equation*}
S(T):=\frac{1}{\pi }\mathrm{Im}\log \zeta (\frac{1}{2}+iT),
\end{equation*}%
where the logarithm is calculated by continuous variation along the contour $%
\sigma +iT,$ with $\sigma $ changing from $+\infty $ to $1/2.$

It is known (see Chapter 15 in \cite{davenport67}) that 
\begin{equation*}
\mathcal{N}(T)=\frac{T}{2\pi }\log \frac{T}{2\pi e}+\frac{7}{8}+S\left(
T\right) +O\left( \frac{1}{1+T}\right) .
\end{equation*}%
Let $t_{k}$ be the solution of the equation 
\begin{equation*}
\frac{t}{2\pi }\log \frac{t}{2\pi e}+\frac{7}{8}=k-1/2.
\end{equation*}%
It is convenient to think about $t_{k}$ as predicted imaginary parts of
Riemann's zeros, $\gamma _{k}.$ Note that the distance between consecutive $%
t_{k}$ are of order $1/\log t_{k}.$ Let 
\begin{equation*}
\sigma _{k}:=\frac{\sqrt{2\log \log t_{k}}}{\log t_{k}},
\end{equation*}%
and define 
\begin{equation}
f_{k}=\frac{\gamma _{k}-t_{k}}{\sigma _{k}}.  \label{definition_f_k}
\end{equation}

The quantities $f_{k}$ show normalized fluctuations of imaginary parts of
Riemann's zeros from their predicted locations $t_{k}$. In order to study
the statistical properties of $f_{k}$ we introduce a probability space $%
\left\{ \Omega ,\mathcal{B},\mathbb{P}\right\} $, where $\Omega =\left[ 0,1%
\right] ,$ $\mathcal{B}$ is the $\sigma $-algebra of Borel subsets of $%
\Omega ,$ and $\mathbb{P}$ is the Lebesgue measure on $\mathcal{B}$.

Let us fix $\theta \in (1/2,1].$ We define a sequence of random variables $%
f^{\left( N\right) }$ by the following formula: 
\begin{equation}
f^{\left( N\right) }\left( \omega \right) :=f_{k\left( N,\omega \right) },
\label{definition_rv}
\end{equation}%
where 
\begin{equation}
k\left( N,\omega \right) :=\left\lfloor N+\omega \left\lfloor N^{\theta
}\right\rfloor \right\rfloor ,  \label{definition_k}
\end{equation}%
$\omega \in \Omega ,$ and $\left\lfloor x\right\rfloor $ denotes the largest
integer which is less than or equal to $x.$ Hence $k\left( N,\omega \right) $
is a random variable uniformly distributed on $I_{N}=\mathbb{Z}\cap \left[
N,N+\left\lfloor N^{\theta }\right\rfloor -1\right] .$ Note that 
\begin{equation*}
\mathbb{P}\left\{ f^{\left( N\right) }\in \left( a,b\right) \right\} =\frac{1%
}{\left\lfloor N^{\theta }\right\rfloor }\left| \left\{ k:k\in
I_{N},f_{k}\in \left( a,b\right) \right\} \right| ,
\end{equation*}%
and 
\begin{equation*}
\mathbb{E}\left( f^{\left( N\right) }\right) ^{r}=\frac{1}{\left\lfloor
N^{\theta }\right\rfloor }\sum_{k\in I_{N}}\left( f_{k}\right) ^{r}.
\end{equation*}

First, we will prove the following theorem.

\begin{theorem}
\label{theorem_main}Suppose that random variables $f^{\left( N\right) }$ are
defined as in (\ref{definition_rv}) with $1/2<\theta \leq 1.$ Then, as $%
N\rightarrow \infty ,$ we have: \newline
(i) for every real $\xi ,$ 
\begin{equation*}
\mathbb{P}\left\{ f^{\left( N\right) }>\xi \right\} \rightarrow \frac{1}{%
\sqrt{2\pi }}\int_{\xi }^{\infty }e^{-x^{2}/2}dx,\text{ and}
\end{equation*}%
(ii) for every integer $p\geq 0,$\newline
\begin{equation*}
\lim_{N\rightarrow \infty }\mathbb{E}\left( f^{\left( N\right) }\right) ^{p}=%
\frac{1}{\sqrt{2\pi }}\int_{-\infty }^{\infty }x^{p}e^{-x^{2}/2}dx.
\end{equation*}
\end{theorem}

The requirement that $\theta >1/2$ comes from a density estimate for Riemann
zeros. This estimate says that the number of zeros with the imaginary part
in the interval $\left[ T,T+H\right] $ and the real part in $\left( \sigma
,\infty \right) $, $\sigma \geq 1/2$, is bounded by a multiple of $%
HT^{-\alpha \left( \sigma -1/2\right) }\log T,$ where $\alpha $ is a
positive constant, provided that $T$ is sufficiently large. The bound is
uniform in $\sigma .$ Selberg's density theorem (Theorem 1 in \cite%
{selberg46a}) establishes this result for $H\geq T^{\theta },$ $\theta >1/2.$
Karatsuba (\cite{karatsuba96}, \cite{karatsuba_korolev06}) established the
density estimate for $\theta >27/82.$ Moreover, Korolev showed in \cite%
{korolev03} that the density estimate holds for ``almost all'' $T$ if $%
H>T^{\varepsilon },$ where $\varepsilon $ is an arbitrary positive constant.
We expect that the results in our Theorem \ref{theorem_main} can be improved
to include the cases $\theta >27/82$ and perhaps even the case $\theta >0$
by using these density estimates.

It is also interesting to ask how $f_{k}$ and $f_{k^{\prime }}$ are related
when $k$ and $k^{\prime }$ are sufficiently close to each other. More
precisely, define random variables $f_{1}^{\left( N\right) }$ and $%
f_{2}^{\left( N\right) }$ by the formula: 
\begin{equation}
f_{i}^{\left( N\right) }\left( \omega \right) =f_{k_{i}\left( N,\omega
\right) },\text{ }i=1,2,  \label{definition_rv_f1_f2}
\end{equation}%
where%
\begin{equation*}
k_{1}\left( N,\omega \right) =N+\left\lfloor \omega N\right\rfloor ,
\end{equation*}%
\begin{equation*}
k_{2}\left( N,\omega \right) =N+\left\lfloor \omega N\right\rfloor +\left[
\left( \log N\right) ^{\beta }\right] ,
\end{equation*}%
and $\beta >0$. (We set here $\theta =1$ for simplicity. However, the result
below is likely to hold for all $\theta \in (27/82,1].$)

\begin{theorem}
\label{theorem_main2}Suppose that random variables $f_{i}^{\left( N\right) }$
are defined as in (\ref{definition_rv_f1_f2}) with $\beta >0,$ and suppose
that $Y_{1},Y_{2}$ are zero-mean Gaussian random variables with $\mathbb{E}%
\left( Y_{i}^{2}\right) =1$ and $\mathbb{E}\left( Y_{1}Y_{2}\right) =\left(
1-\beta \right) _{+}.$ Then as $N\rightarrow \infty $, \newline
(i) the joint cumulative distribution function of $\left( f_{1}^{\left(
N\right) },f_{2}^{\left( N\right) }\right) $ converges pointwise to the
joint cumulative distribution function of $\left( Y_{1},Y_{2}\right) $, and 
\newline
(ii) the joint moments of $\left( f_{1}^{\left( N\right) },f_{2}^{\left(
N\right) }\right) $ converge to the corresponding joint moments of $\left(
Y_{1},Y_{2}\right) .$
\end{theorem}

After the first version of this article was completed, the author learned
from M. A. Korolev about his papers \cite{korolev10a} and \cite{korolev10b}
(based on earlier results by Karatsuba and Korolev in \cite{karatsuba96} and %
\cite{karatsuba_korolev06}), that consider similar questions. See, for
example, Theorem 10 in \cite{korolev10b} which is similar to our Theorem \ref%
{theorem_main}. However, the joint distribution of the fluctuations of zeta zeros is
not studied in these papers.

We have shown that the distribution of two zeta zero fluctuations approaches
a two-variate Gaussian distribution. By a natural extension of the argument,
with a more cumbersome notation, it is possible to show that the
distribution of any finite number of fluctuations approaches a multivariate
Gaussian distribution with the covariance matrix $\mathbb{E}%
X_{i}X_{j}=\left( 1-\beta _{ij}\right) _{+},$ where 
\begin{equation*}
\beta _{ij}=\lim_{N\rightarrow \infty }\frac{\log \left| k_{j}\left(
N,\omega \right) -k_{i}\left( N,\omega \right) \right| }{\log \log N},
\end{equation*}%
and the limit is assumed to be positive and the same for all $\omega .$ (The
fluctuations are around the predicted locations $t_{k_{i}(N,\omega )},$ and
the functions $k_{i}\left( N,\omega \right) $ are defined as in (\ref%
{definition_rv_f1_f2}) with appropriate changes.)

The positive numbers $\beta _{ij}$ in the covariance matrix are not
arbitrary but satisfy the ultrametric inequality: 
\begin{equation*}
\beta _{ik}\leq \max \left\{ \beta _{ij},\beta _{jk}\right\} .
\end{equation*}%
More about this covariance structure can be found in Section 4 of \cite%
{bourgade10}, where it is shown, in particular, how this structure can arise
as a result of a branching process.

Covariances that satisfy ultrametric inequalities are of interest in
statistical physics. They are used, in particular, in the theory of
frustrated disordered systems (``spin glasses''), where they are crucial in
a proposed description of local equilibria by replica method (see \cite%
{rtv86}, \cite{mpv87}, and \cite{talagrand11}). A possible reason for the
appearance of ultrametric structure in this area of physics is the close
relation of spin glasses with random matrices where ultrametric covariances
describe the eigenvalue distribution.

In particular, Diaconis and Evans in \cite{diaconis_evans01} considered
uniformly distributed $N$-by-$N$ random unitary matrices and determined the
covariances for the eigenvalue counts in given intervals for large $N$ (see
Theorems 6.1 and 6.3 in their paper). They found that these covariances have
an interesting and unusual structure. A similar structure was found for
Gaussian Hermitian random matrices and other random matrix ensembles by
Soshnikov in \cite{soshnikov00a}. In fact, for random unitary matrices this
structure can be seen as a consequence of the ultrametric covariances
exhibited by characteristic polynomials of these matrices (Theorem
1.4 in Bourgade's paper \cite{bourgade10}). Bourgade has also found a
parallel result for counts of Riemann zeros in given intervals (Theorem 1.1
and Corollary 1.3 in \cite{bourgade10}).

In another development, Gustavsson (\cite{gustavsson05} ) studied
eigenvalues of Gaussian Hermitian matrices and found the ultrametric
structure in covariances defined by using the deviations of individual
eigenvalues from their predicted locations (Theorems 1.3 and 1.4 in \cite%
{gustavsson05}). This setup is similar to what we do in this paper and the
results are also remarkably similar.

However, while the results are similar, the methods are quite different. In
random matrix papers, the method is based either on group representation
theory which allows one to compute average traces of matrix powers (as in %
\cite{diaconis_evans01}, and \cite{bourgade10}), or on explicit formulas for
the distribution of eigenvalues (as in \cite{soshnikov00a}, \cite%
{gustavsson05}). In contrast, in number-theoretic papers, the method is
based on the Selberg approximation formula for the number of Riemann zeros
with ordinates between zero and T, and on a multitude of other facts from
number theory, which allow one to estimate the powers of this approximate
function. The fact that these distinct methods lead to very similar results
is rather mysterious.

An interested reader can find more about relations of Riemann's zeros and
random matrices in review papers mentioned in the beginning of this paper.

The rest of the paper is organized as follows. Section \ref%
{section_outline_of_proofs} outlines the scheme of the proof of Theorems \ref%
{theorem_main} and \ref{theorem_main2}. Section \ref{section_exp_sums}
introduces some technical tools that we will need in the proof of the main
theorems. Section \ref{section_Selberg_approximation} proves a modification
of the key approximation result by Selberg. Section \ref{section_moments_Sx}
calculates the moments of the approximate function $S_{x}.$ Section \ref%
{section_moments_S_and_proof_of_Thm1} calculates the moments of $S$ and
concludes the proof of Theorem \ref{theorem_main}. Section \ \ref%
{section_covariance} proves Theorem \ref{theorem_covariance}. Section \ref%
{section_joint_moments} proves Theorem \ref{theorem_joint_moments_X} and
concludes the proof of Theorem \ref{theorem_main2}. And Section \ref%
{section_conclusion} concludes.

\section{Outline of Proofs}

\label{section_outline_of_proofs}

In the proof we use the strategy used by Gustavsson in his work on the
fluctuations of eigenvalues in the Gaussian Unitary Ensemble. The first step
in Gustavsson's proof is to relate fluctuations of an individual eigenvalue
to fluctuations of eigenvalue counts in a fixed interval. This allows one to
use existing methods for finding the distribution of eigenvalue counts.

In our setup, an analogous step requires connecting the random fluctuations $%
f_{k\left( N,\omega \right) }$ to the number of Riemann zeros in the
interval $\left[ 0,t_{k\left( N,\omega \right) }\right] ,$ which can be
approximated by the function $S\left( t_{k\left( N,\omega \right) }\right) .$

It is convenient to define 
\begin{equation*}
X_{k}:=\frac{\sqrt{2}\pi S\left( t_{k}+\xi \sigma _{k}\right) }{\sqrt{\log
\log t_{k}}},
\end{equation*}%
and a corresponding sequence of random variables 
\begin{equation}
X^{\left( N\right) }\left( \omega \right) :=X_{k\left( N,\omega \right) },
\label{definition_rv_X}
\end{equation}%
where $k\left( N,\omega \right) $ is as in (\ref{definition_k}).

A connection between $X^{(N)}$ and $f^{(N)}$ can be seen as follows. For
every real $\xi ,$ 
\begin{eqnarray*}
\mathbb{P}\left\{ f^{\left( N\right) }>\xi \right\} &=&\frac{1}{H_{N}}%
\left\vert \left\{ k:k\in I_{N},\gamma _{k}>t_{k}+\xi \sigma _{k}\right\}
\right\vert \\
&=&\frac{1}{H_{N}}\left\vert \left\{ k:k\in I_{N},\mathcal{N}\left(
t_{k}+\xi \sigma _{k}\right) \leq k-1/2\right\} \right\vert \\
&=&\frac{1}{H_{N}}\left\vert \left\{ k:k\in I_{N},\frac{t_{k}+\xi \sigma _{k}%
}{2\pi }\log \frac{t_{k}+\xi \sigma _{k}}{2\pi e}+\frac{7}{8}+S\left(
t_{k}+\xi \sigma _{k}\right) +O\left( 1/t_{k}\right) \leq k-1/2\right\}
\right\vert \\
&=&\frac{1}{H_{N}}\left\vert \left\{ k:k\in I_{N},S\left( t_{k}+\xi \sigma
_{k}\right) \leq -\xi \sqrt{\frac{\log \log t_{k}}{2\pi ^{2}}}\left( 1+\frac{%
\log 2\pi }{\log t_{k}}\right) +o\left( 1/t_{k}\right) \right\} \right\vert
\\
&=&\frac{1}{H_{N}}\left\vert \left\{ k:k\in I_{N},\frac{\sqrt{2}\pi S\left(
t_{k}+\xi \sigma _{k}\right) }{\sqrt{\log \log t_{k}}}\leq -\xi \left( 1+%
\frac{\log 2\pi }{\log t_{k}}\right) +o\left( 1/t_{k}\right) \right\}
\right\vert .
\end{eqnarray*}

Since $t_{k}$ is asymptotically close to $2\pi k/\log k,$ it follows that 
\begin{equation}
\mathbb{P}\left\{ X^{\left( N\right) }\leq -\xi -\frac{c_{1}\left| \xi
\right| }{\log N}+o\left( \frac{\log N}{N}\right) \right\} \leq \mathbb{P}%
\left\{ f^{\left( N\right) }>\xi \right\} \leq \mathbb{P}\left\{ X^{\left(
N\right) }\leq -\xi +\frac{c_{2}\left| \xi \right| }{\log N}+o\left( \frac{%
\log N}{N}\right) \right\} ,  \label{bounds_distr_f}
\end{equation}%
where $c_{1}$ and $c_{2}$ are two constants. Hence for large $N,$ the
distribution of the random variable $f^{\left( N\right) }$ is essentially
determined by the distribution of the random variable $X^{\left( N\right) }.$

In this connection, it is appropriate to recall the following theorem by
Selberg (Theorem 3 in \cite{selberg44}). Let 
\begin{equation*}
X\left( t\right) :=\frac{\sqrt{2}\pi S(t)}{\sqrt{\log \log t}}.
\end{equation*}

\begin{theorem}[Selberg]
\label{theorem_selberg} Assume RH, and let $T^{a}\leq H\leq T^{2},$ where $%
a>0.$ Then for every $k\geq 1$ \ 
\begin{equation*}
\frac{1}{H}\int_{T}^{T+H}\left| X\left( t\right) \right| ^{2k}dt=\frac{2k!}{%
k!2^{k}}+O(1/\log \log T),
\end{equation*}%
with the constant in the remainder term that depends only on $k$ and $a.$
\end{theorem}

In other words, the even moments of the function $X\left( t\right) $ behave
as the moments of a standard Gaussian variable. This was refined in \cite%
{ghosh83}, where it was shown in particular that for every interval $I,$%
\begin{equation*}
\frac{1}{H}\int_{T}^{T+H}\mathbf{1}_{I}\left[ X\left( t\right) \right] dt=%
\frac{1}{\sqrt{2\pi }}\int_{I}e^{-x^{2}/2}dx+o(1),
\end{equation*}%
where $\mathbf{1}_{I}$ denotes the indicator function of interval $I.$ We
prove a modified version of this result.

\begin{theorem}
\label{theorem_main_X}Suppose that random variables $X^{\left( N\right) }$
are defined as in (\ref{definition_rv_X}) with $1/2<\theta \leq 1.$ Then,
for every real $s,$ as $N\rightarrow \infty ,$%
\begin{equation*}
\mathbb{E}\mathbf{1}_{\left( -\infty ,s\right] }\left( X^{\left( N\right)
}\right) \rightarrow \frac{1}{\sqrt{2\pi }}\int_{-\infty }^{s}e^{-x^{2}/2}dx.
\end{equation*}
\end{theorem}

We will prove Theorem \ref{theorem_main_X} by the method of moments, which
says that in order to establish the convergence of a sequence of r.v. in
distribution to the Gaussian law it is enough to show the convergence of
every moment (Example 2.23 on p. 18 in van der Vaart \cite{van_der_vaart98}%
). That is, it is enough to show that 
\begin{equation}
\mathbb{E}\left( X^{\left( N\right) }\right) ^{r}\rightarrow \frac{1}{\sqrt{%
2\pi }}\int_{-\infty }^{\infty }x^{r}e^{-x^{2}/2}dx
\label{convergence_moments_XN}
\end{equation}%
for every integer $r>0.$ We will show this in Section \ref%
{section_moments_S_and_proof_of_Thm1} in Corollary \ref{corollary_moments}.

The first claim in Theorem \ref{theorem_main} follows \ immediately from
Theorem \ref{theorem_main_X} and inequalities (\ref{bounds_distr_f}). The
second claim follows from the first one because (\ref{bounds_distr_f}) and (%
\ref{convergence_moments_XN}) imply that $\left( f^{\left( N\right) }\right)
^{2n}$ are asymptotically uniformly integrable for every $n>0$ and therefore
the moments of $f^{\left( N\right) }$ converge to the moments of the
limiting Gaussian distribution (see Theorem 2.20 in van der Vaart \cite%
{van_der_vaart98}).

In order to prove Theorem \ref{theorem_main2}, define random variables 
\begin{equation}
X_{i}^{\left( N\right) }\left( \omega \right) :=X_{k_{i}\left( N,\omega
\right) },\text{ }i=1,2,  \label{definition_rv_Xi}
\end{equation}%
where $k_{i}\left( N,\omega \right) $ are as in (\ref{definition_rv_f1_f2}).
Let us use notations \newline
$I_{1}^{\left( N\right) }:=[N,2N-1],$ $I_{2}^{\left( N\right) }:=\left[ N+%
\left[ \left( \log N\right) ^{\beta }\right] ,2N-1+\left[ \left( \log
N\right) ^{\beta }\right] \right] ,$ and $\mathbf{s}:=\left(
s_{1},s_{2}\right) .$ Then, 
\begin{eqnarray*}
\mathbb{P}\left\{ f_{1}^{\left( N\right) }>\xi _{1},f_{2}^{\left( N\right)
}>\xi _{2}\right\} &=&N^{-2}\left| \left\{ \mathbf{s\in }\mathbb{Z}%
^{2}:s_{i}\in I_{i}^{\left( N\right) },\gamma _{s_{i}}>t_{s_{i}}+\xi
_{i}\sigma _{s_{i}},\text{ }i=1,2\right\} \right| \\
&=&N^{-2}\left| \left\{ \mathbf{s}:s_{i}\in I_{i}^{\left( N\right) },%
\mathcal{N}\left( t_{s_{i}}+\xi _{i}\sigma _{s_{i}}\right) \leq
k_{i}-1/2\right\} \right| \\
&=&N^{-2}\left| \left\{ \mathbf{s}:s_{i}\in I_{i}^{\left( N\right) },\frac{%
\sqrt{2}\pi S\left( t_{s_{i}}+\xi _{i}\sigma _{s_{i}}\right) }{\sqrt{\log
\log t_{s_{i}}}}\leq -\xi _{i}\left( 1+\frac{\log 2\pi }{\log t_{s_{i}}}%
\right) +o\left( 1/t_{s_{i}}\right) \right\} \right| .
\end{eqnarray*}

That is, with some positive $c_{1}$ and $c_{2},$ we have 
\begin{equation}
\mathbb{P}\left\{ f_{1}^{\left( N\right) }>\xi _{1},f_{2}^{\left( N\right)
}>\xi _{2}\right\} \leq \mathbb{P}\left\{ X_{i}^{\left( N\right) }\leq -\xi
_{i}+\frac{c_{1}}{\log N}\left| \xi _{i}\right| +o\left( \frac{\log N}{N}%
\right) ,\text{ }i=1,2\right\} ,  \label{multivariate_bound1}
\end{equation}%
and 
\begin{equation}
\mathbb{P}\left\{ f_{1}^{\left( N\right) }>\xi _{1},f_{2}^{\left( N\right)
}>\xi _{2}\right\} \geq \mathbb{P}\left\{ X_{i}^{\left( N\right) }\leq -\xi
_{i}-\frac{c_{2}}{\log N}\left| \xi _{i}\right| +o\left( \frac{\log N}{N}%
\right) ,\text{ }i=1,2\right\} .  \label{multivariate_bound2}
\end{equation}

In words, the joint cumulative distribution function of $f_{1}^{\left(
N\right) }$ and $f_{2}^{\left( N\right) }$ approaches that of $X_{1}^{\left(
N\right) }$ and $X_{2}^{\left( N\right) }$.

First of all, we have the following result for the random variables $%
X_{1}^{\left( N\right) }$ and $X_{2}^{\left( N\right) }.$

\begin{theorem}
\label{theorem_covariance} Let $X_{i}^{\left( N\right) }$ be defined as in (%
\ref{definition_rv_Xi}). Then, 
\begin{equation*}
\lim_{N\rightarrow \infty }\mathbb{E}X_{1}^{\left( N\right) }X_{2}^{\left(
N\right) }=\left( 1-\beta \right) _{+}:=\left\{ 
\begin{array}{cc}
1-\beta , & \text{ if }\beta \in \left( 0,1\right) , \\ 
0, & \text{if }\beta \geq 1.%
\end{array}%
\right.
\end{equation*}
\end{theorem}

More generally, the following result holds.

\begin{theorem}
\label{theorem_joint_moments_X}Let $X_{i}^{\left( N\right) }$ be defined as
in (\ref{definition_rv_Xi}). Then for every $l,m\geq 0,$ $\mathbb{E}\left(
X_{1}^{\left( N\right) }\right) ^{l}\left( X_{2}^{\left( N\right) }\right)
^{m}$ converges to $\mathbb{E}\left( Y_{1}\right) ^{l}\left( Y_{2}\right)
^{m}$ where $\left( Y_{1},Y_{2}\right) $ is a zero-mean Gaussian random
variable with $\mathbb{E}\left( Y_{i}^{2}\right) =1$ and $\mathbb{E}\left(
Y_{1}Y_{2}\right) =\left( 1-\beta \right) _{+}.$
\end{theorem}

Theorem \ref{theorem_covariance} is a particular case of Theorem \ref%
{theorem_joint_moments_X}. However, we will prove it separately, since its
proof is more transparent and shows how the proof of the more general
Theorem \ref{theorem_joint_moments_X} proceeds.

Given Theorem \ref{theorem_joint_moments_X}, we can prove Theorem \textbf{\ }%
\ref{theorem_main2}.

\textbf{Proof of Theorem \ref{theorem_main2}:} Theorem \ref%
{theorem_joint_moments_X} implies that the cumulative distribution function
of $\left( X_{1}^{\left( N\right) },X_{2}^{\left( N\right) }\right) $
converges pointwise to the cumulative distribution function of the Gaussian
variable $\left( X_{1},X_{2}\right) $. The first claim of the theorem
follows immediately from this fact and inequalities (\ref%
{multivariate_bound1}) and (\ref{multivariate_bound2}). In addition, Theorem %
\ref{theorem_joint_moments_X} and inequalities (\ref{multivariate_bound1})
and (\ref{multivariate_bound2}) imply that for all integer $a,b\geq 0,$ the
random variables $\left( f_{1}^{\left( N\right) }\right) ^{a}\left(
f_{2}^{(N)}\right) ^{b}$ are asymptotically uniformly integrable. Hence,
their expectations converge to the corresponding expectation of the limit, $%
\mathbb{E}\left( Y_{1}\right) ^{a}\left( Y_{2}\right) ^{b}$ (by Theorem 2.20
in van der Vaart \cite{van_der_vaart98}). This completes the proof of the
second claim of the theorem. $\square $

The proof of the convergence of moments of $X^{\left( N\right) }$ follows
the plan of the argument in Selberg \cite{selberg44}.

Recall that $X^{\left( N\right) }$ is a rescaled version of $S\left(
t_{k}+\xi \sigma _{k}\right) $ where $k$ is random. The first step in
Selberg's proof is to show that $S\left( t\right) $ can be approximated by $%
S_{x}\left( t\right) ,$ where 
\begin{equation*}
S_{x}\left( t\right) :=-\frac{1}{\pi }\sum_{p\leq x^{3}}\frac{\sin \left(
t\log p\right) }{\sqrt{p}}.
\end{equation*}%
That is, Selberg shows that 
\begin{equation*}
\frac{1}{H}\int_{K}^{K+H}\left( S\left( t\right) -S_{x}\left( t\right)
\right) ^{2n}dt
\end{equation*}%
is small provided that $K$ and $H$ are sufficiently large and that $x\sim
K^{\varepsilon }$ with a sufficiently small $\varepsilon >0.$ In our case we
will need to modify this result in order to show that the integral can be
replaced by a sum over a discrete set of points.

The next step in Selberg's proof is to calculate the moments 
\begin{equation*}
\frac{1}{H}\int_{K}^{K+H}\left\vert S_{x}\left( t\right) \right\vert ^{2n}dt.
\end{equation*}%
Again it will be necessary to prove a corresponding result for a sum over a
discrete set of points.

Given the results in these two steps, it is relatively easy to calculate the
moments of the random variable $S\left( t_{k}+\xi \sigma _{k}\right) .$ This
will be done essentially as in Selberg's paper. However, we will need to
extend the calculation to the multivariate case with two random variables $%
S\left( t_{k}+\xi _{1}\sigma _{k}\right) $ and $S\left( t_{k}+\xi _{2}\sigma
_{k}\right) $

\section{Exponential Sums}

\label{section_exp_sums}

The changes in Selberg's proof make it necessary to estimate certain
exponential sums. The main additional tool that we use to handle these sums
is the following theorem by van der Corput (Theorem 2.2 in \cite%
{graham_kolesnik91}). Let $e\left( f\left( n\right) \right) $ denote $\exp %
\left[ 2\pi if\left( n\right) \right] .$

\begin{theorem}[van der Corput]
\label{theorem_van_der_Corput}Suppose that $f$ is a real valued function
with two continuous derivatives on interval $I.$ Suppose also that there is
some $\lambda >0$ and some $\kappa \geq 1$ such that 
\begin{equation*}
\lambda \leq \left| f^{\prime \prime }\left( x\right) \right| \leq \kappa
\lambda
\end{equation*}%
on $I.$ Then,%
\begin{equation*}
\sum_{n\in I}e\left( f\left( n\right) \right) =O\left( \kappa \left|
I\right| \lambda ^{1/2}+\lambda ^{-1/2}\right) .
\end{equation*}
\end{theorem}

In order to apply this theorem in our situation, we need to estimate
derivatives of a function $g(x)$ that we are about to define. Let $t\left(
x\right) $ be the functional inverse of the function 
\begin{equation}
x(t)=\frac{t}{2\pi }\log \frac{t}{2\pi e}+\frac{11}{8}  \label{formula_x}
\end{equation}%
on interval $[t_{0},\infty )$ where $t_{0}$ is sufficiently large. (The
coefficients $11/8$ is not necessary is not necessary for the argument. It
is included only because its presence makes $t(k)$ an unbiased estimator of $%
\gamma _{k}.$) Note that $t\left( x\right) $ is an increasing concave
function. Let 
\begin{equation}
g\left( x\right) :=t\left( x\right) +\xi \frac{\sqrt{2\log \log t\left(
x\right) }}{\log t\left( x\right) },  \label{definition_g}
\end{equation}%
where $\xi $ is a real constant. This function is well defined for $x$
greater than some numeric constant $x_{0}.$ For $x$ between $0$ and $x_{0}$,
we define $g(x)$ in an arbitrary fashion such that $g(x)$ has a continuous $%
3 $-rd derivative for all $x\geq 0$.

\begin{lemma}
\label{lemma_derivatives_g}When $x\rightarrow \infty $, 
\begin{equation*}
g^{\prime \prime }(x)\sim -\frac{2\pi }{x\left( \log x\right) ^{2}}
\end{equation*}%
and 
\begin{equation*}
g^{\prime \prime \prime }\left( x\right) \sim \frac{2\pi }{x^{2}\left( \log
x\right) ^{2}}.
\end{equation*}
\end{lemma}

\textbf{Proof:} The identity 
\begin{equation*}
x=\frac{t\left( x\right) }{2\pi }\log \frac{t\left( x\right) }{2\pi e}+\frac{%
11}{8}
\end{equation*}%
implies that

\begin{equation*}
t\sim 2\pi \frac{x}{\log x}\text{, }t^{\prime }\sim 2\pi \frac{1}{\log x},%
\text{ }t^{\prime \prime }\sim -\frac{2\pi }{x\left( \log x\right) ^{2}},%
\text{ and }t^{\prime \prime \prime }\sim \frac{2\pi }{x^{2}\left( \log
x\right) ^{2}}.
\end{equation*}

If 
\begin{equation*}
h:=\frac{\sqrt{\log \log t\left( x\right) }}{\log t\left( x\right) },
\end{equation*}%
then a calculation shows that $h^{\prime \prime }=o(t^{\prime \prime }),$ $%
h^{\prime \prime \prime }=o(t^{\prime \prime \prime })$ and therefore 
\begin{equation*}
g^{\prime \prime }\sim -\frac{2\pi }{x\left( \log x\right) ^{2}}\text{ and }%
g^{\prime \prime \prime }\sim \frac{2\pi }{x^{2}\left( \log x\right) ^{2}}.
\end{equation*}%
$\square $

In the following we will use the notation $g_{k}$ for $g\left( k\right)
\equiv t_{k}+\xi \sigma _{k}.$

\begin{lemma}
\label{lemma_estimate_exp_sum}Let 
\begin{equation*}
\theta =\log \left( \frac{p_{l+1}\ldots p_{2n}}{p_{1}\ldots p_{l}}\right) ,
\end{equation*}%
where $1\leq l\leq 2n,$ $\left\{ p_{1},\ldots ,p_{l}\right\} \neq \left\{
p_{l+1},\ldots p_{2n}\right\} $ and primes $p_{i}<y$ for all $i$. Assume $%
1\leq H\leq cK.$ Then, 
\begin{equation*}
\sum_{k=K}^{K+H-1}e^{i\theta g_{k}}=O\left( Hn\frac{y^{n/2}\log y}{%
K^{1/2}\log K}+y^{n/2}K^{1/2}\log K\right) .
\end{equation*}
\end{lemma}

\textbf{Proof}: From the assumption, we obtain%
\begin{equation*}
c/y^{n}\leq \left\vert \theta \right\vert \leq 2n\log y.
\end{equation*}%
(In order to see the first inequality, let $l\leq n.$ Then 
\begin{equation*}
\left\vert 1-\frac{p_{l+1}\ldots p_{2n}}{p_{1}\ldots p_{l}}\right\vert
=\left\vert \frac{p_{1}\ldots p_{l}-p_{l+1}\ldots p_{2n}}{p_{1}\ldots p_{l}}%
\right\vert \geq \frac{1}{y^{n}}
\end{equation*}%
by the uniqueness of integer factorization, and the desired inequality
follows. The case $l\geq n$ is similar.)

Hence, by using Lemma \ref{lemma_derivatives_g}, we find that 
\begin{equation*}
\lambda \leq \left\vert \theta g^{\prime \prime }\left( x\right) \right\vert
\leq \kappa \lambda
\end{equation*}%
with 
\begin{equation*}
\lambda =\frac{c}{y^{n}K\left( \log K\right) ^{2}},
\end{equation*}%
and 
\begin{equation*}
\kappa =O(ny^{n}\log y)
\end{equation*}%
By applying van der Corput's theorem, we obtain 
\begin{equation*}
\sum_{k=K}^{K+H-1}e^{i\theta g_{k}}=O\left( Hn\frac{y^{n/2}\log y}{%
K^{1/2}\log K}+y^{n/2}K^{1/2}\log K\right) .
\end{equation*}%
$\square $

\begin{lemma}
\label{lemma_bound_p_05}Suppose $1\leq c_{1}K^{\theta }\leq H\leq c_{2}K,$
where $\theta >1/2$ and $c_{1},c_{2}>0.$ Let $r$ be a positive integer, $%
y\leq K^{\frac{2\theta -1}{3r}-\varepsilon },$ and assume that 
\begin{equation*}
\left| \alpha _{p}\right| <A\frac{\log p}{\log y}\text{ for }p<y.
\end{equation*}%
Then, we have 
\begin{equation*}
\sum_{k=K}^{K+H-1}\left| \sum_{p<y}\frac{\alpha _{p}}{p^{1/2+ig_{k}}}\right|
^{2r}=O\left( H\right) .
\end{equation*}
\end{lemma}

\textbf{Proof:} We can write 
\begin{equation*}
\left( \sum_{p<y}\frac{\alpha _{p}}{p^{1/2+ig_{k}}}\right)
^{r}=\sum_{n<y^{r}}\frac{\beta _{n}}{n^{1/2+ig_{k}}},
\end{equation*}%
where $\beta _{n}\leq A^{r}$. Hence, 
\begin{eqnarray*}
\sum_{k=K}^{K+H-1}\left| \sum_{p<y}\frac{\alpha _{p}}{p^{1/2+ig_{k}}}\right|
^{2r} &=&\sum_{m,n<y^{r}}\frac{\beta _{m}\overline{\beta }_{n}}{\sqrt{mn}}%
\sum_{k=K}^{K+H-1}\left( \frac{m}{n}\right) ^{ig_{k}} \\
&\leq &H\sum_{n<y^{r}}\frac{\left| \beta _{n}\right| ^{2}}{n}%
+2\sum_{m<n<y^{r}}\frac{\left| \beta _{m}\beta _{n}\right| }{\sqrt{mn}}%
\left| \sum_{k=K}^{K+H-1}\left( \frac{m}{n}\right) ^{ig_{k}}\right| .
\end{eqnarray*}%
The first sum can be estimated as follows: 
\begin{equation*}
\sum_{n<y^{r}}\frac{\left| \beta _{n}\right| ^{2}}{n}\leq A^{r}\sum_{n<y^{r}}%
\frac{\left| \beta _{n}\right| }{n}\leq A^{r}\left( \sum_{p<y}\frac{\left|
\alpha _{p}\right| }{p}\right) ^{r}=O(1),
\end{equation*}%
where we used Mertens' result $\sum_{p<y}\frac{\log p}{p}=O(\log y)$ in the
last step.

In order to estimate the second sum we note that 
\begin{equation*}
1/y^{r}<\log \left\vert n/m\right\vert <r\log y;
\end{equation*}%
hence we can apply van der Corput's theorem and estimate 
\begin{equation*}
\left\vert \sum_{k=K}^{K+H-1}\left( \frac{m}{n}\right) ^{ig_{k}}\right\vert
\leq O\left( Hr\frac{y^{r/2}\log y}{K^{1/2}\log K}+y^{r/2}K^{1/2}\log
K\right) .
\end{equation*}%
Besides, 
\begin{equation*}
\sum_{m<n<y^{r}}\frac{\left\vert \beta _{m}\beta _{n}\right\vert }{\sqrt{mn}}%
\leq \left( \sum_{p<y}\frac{\alpha _{p}}{\sqrt{p}}\right) ^{2r}=O\left(
y^{r}\right) .
\end{equation*}

By assumptions about $H$ and $y,$ it follows that%
\begin{eqnarray*}
\sum_{m<n<y^{r}}\frac{\left\vert \beta _{m}\beta _{n}\right\vert }{\sqrt{mn}}%
\left\vert \sum_{k=K}^{K+H-1}\left( \frac{m}{n}\right) ^{ig_{k}}\right\vert
&=&O\left( y^{r}\right) O\left( Hr\frac{y^{r/2}\log y}{K^{1/2}\log K}%
+y^{r/2}K^{1/2}\log K\right) \\
&=&O(H).
\end{eqnarray*}
$\square $

\begin{lemma}
\label{lemma_bound_p_1}Suppose $1\leq c_{1}K^{\theta }\leq H\leq c_{2}K,$
where $\theta >1/2$ and $c_{1},c_{2}>0.$ Let $r$ be a positive integer, $%
y\leq K^{\frac{2\theta -1}{3r}-\varepsilon },$ and assume that 
\begin{equation*}
\left| \alpha _{p}\right| <A\text{ for }p<y.
\end{equation*}%
Then, 
\begin{equation*}
\sum_{k=K}^{K+H-1}\left| \sum_{p<y}\frac{\alpha _{p}}{p^{1+ig_{k}}}\right|
^{2r}=O\left( H\right) .
\end{equation*}
\end{lemma}

The\textbf{\ }proof of this lemma is similar to the proof of the previous
one.

\section{A consequence of Selberg's approximation formula}

\label{section_Selberg_approximation}

Recall that 
\begin{equation*}
S_{x}\left( t\right) :=-\frac{1}{\pi }\sum_{p\leq x^{3}}\frac{\sin \left(
t\log p\right) }{\sqrt{p}}.
\end{equation*}%
Our goal in this section is to prove the following result.

\begin{proposition}
\label{proposition_approximation}Suppose $1\leq c_{1}K^{\theta }\leq H\leq
c_{2}K,$ where $1/2<\theta \leq 1$ and $c_{1},c_{2}>0.$ Let $x=K^{\frac{%
\theta -1/2}{20n}}.$ Then, we have 
\begin{equation*}
\sum_{k=K}^{K+H-1}\left| S\left( g_{k}\right) -S_{x}\left( g_{k}\right)
\right| ^{2n}=O\left( H\right) .
\end{equation*}
\end{proposition}

\textbf{Proof:} Let $\Lambda \left( n\right) =\log p,$ if $n$ is a power of
the prime number $p,$ and $\Lambda \left( n\right) =0,$ otherwise. Also,
define 
\begin{equation*}
\Lambda _{x}\left( n\right) =\left\{ 
\begin{array}{cc}
\Lambda \left( n\right) , & \text{ for }1\leq n\leq x, \\ 
\Lambda \left( n\right) \left( \frac{\log ^{2}\frac{x^{3}}{n}-2\log ^{2}%
\frac{x^{2}}{n}}{2\log ^{2}x}\right) , & \text{for }x\leq n\leq x^{2}, \\ 
\Lambda \left( n\right) \frac{\log ^{2}\frac{x^{3}}{n}}{2\log ^{2}x}, & 
\text{for }x^{2}\leq n\leq x^{3}.%
\end{array}%
\right.
\end{equation*}

Let $a\in (1/2,1],$ $x=T^{\frac{a-1/2}{60k}},$ $T^{a}\leq H\leq T,$ $T\leq
t\leq T+H.$

The first formula on p. 37 in \cite{selberg46a} (immediately before formula
(5.2)) states that 
\begin{eqnarray*}
S\left( t\right) -S_{x}\left( t\right) &=&O\left( \left| \sum_{p<x^{3}}\frac{%
\Lambda \left( p\right) -\Lambda _{x}\left( p\right) }{\sqrt{p}\log p}%
p^{-it}\right| \right) \\
&&+O\left( \left| \sum_{p<x^{3/2}}\frac{\Lambda _{x}\left( p^{2}\right) }{%
p\log p}p^{-2it}\right| \right) +O\left( \left( \sigma _{x,t}-\frac{1}{2}%
\right) \log T\right) \\
&&+O\left( \left( \sigma _{x,t}-\frac{1}{2}\right) x^{\sigma _{x,t}-\frac{1}{%
2}}\int_{1/2}^{\infty }x^{1/2-\sigma }\left| \sum_{p<x^{3}}\frac{\Lambda
_{x}\left( p\right) \log \left( xp\right) }{p^{\sigma +it}}\right| d\sigma
\right) ,
\end{eqnarray*}%
where 
\begin{equation*}
\sigma _{x,t}=\frac{1}{2}+\max_{\rho }\left( \beta -\frac{1}{2},\frac{2}{%
\log x}\right)
\end{equation*}%
and the maximum is taken over all zeros $\beta +i\gamma $ for which 
\begin{equation*}
\left| t-\gamma \right| \leq \frac{x^{3\left| \beta -1/2\right| }}{\log x}.
\end{equation*}

It follows that 
\begin{eqnarray}
&&\sum_{k=K}^{K+H-1}\left| S(g_{k})-S_{x}(g_{k})\right| ^{2n}
\label{Selberg_basic_estimate} \\
&=&O\left( \sum_{k=K}^{K+H-1}\left| \sum_{p<x^{3}}\frac{\Lambda (p)-\Lambda
_{x}(p)}{\sqrt{p}\log p}p^{-ig_{k}}\right| ^{2n}\right)  \notag \\
&&+O\left( \sum_{k=K}^{K+H-1}\left| \sum_{p<x^{3/2}}\frac{\Lambda _{x}(p^{2})%
}{p\log p}p^{-2ig_{k}}\right| ^{2n}\right)  \notag \\
&&+O\left( \left( \log K\right) ^{2n}\sum_{k=K}^{K+H-1}\left( \sigma
_{x,g_{k}}-\frac{1}{2}\right) ^{2n}\right)  \notag \\
&&+O\left( \sum_{k=K}^{K+H-1}\left( \sigma _{x,g_{k}}-\frac{1}{2}\right)
^{2n}x^{2n\left( \sigma _{x,g_{k}}-\frac{1}{2}\right) }\left\{
\int_{1/2}^{\infty }x^{1/2-\sigma }\left| \sum_{p<x^{3}}\frac{\Lambda
_{x}\left( p\right) \log \left( xp\right) }{p^{\sigma +ig_{k}}}\right|
d\sigma \right\} ^{2n}\right) .  \notag
\end{eqnarray}%
By applying Lemmas \ref{lemma_bound_p_05} and \ref{lemma_bound_p_1}, we find
that the first two sums on the right-hand side are $O\left( H\right) .$ For
example, for the first term we can apply Lemma \ref{lemma_bound_p_05} with 
\begin{equation*}
\alpha _{p}=\frac{\Lambda (p)-\Lambda _{x}(p)}{\log p}.
\end{equation*}

The third term can be estimated by using the following lemma.

\begin{lemma}
\label{lemma_bound_integral}Suppose $1\leq c_{1}K^{\theta }\leq H\leq
c_{2}K, $ where $1/2<\theta \leq 1$ and $c_{1},c_{2}>0.$ Next, suppose that $%
x\geq 2, $ $1\leq \xi \leq x^{8k},$ $x^{3}\xi ^{2}\leq \left( \frac{H}{\sqrt{%
K}}\right) ^{1/4}.$ Then we have for $0\leq \nu \leq 8k,$ 
\begin{equation*}
\sum_{k=K}^{K+H-1}\left( \sigma _{x,g_{k}}-\frac{1}{2}\right) ^{\nu }\xi
^{\sigma _{x,g_{k}}-\frac{1}{2}}=O\left( \frac{H}{\left( \log x\right) ^{\nu
}}\right) .
\end{equation*}
\end{lemma}

This lemma is an analog of Lemma 12 on p.33 in \cite{selberg46a} $\ $and its
proof is the same as the proof of Lemma 12 with minor changes. (At this
step, Selberg's density estimate is used.)

By applying Lemma \ref{lemma_bound_integral} with $\xi =1$ and $\nu =2n$ we
find that the third term is 
\begin{equation*}
O\left( H\left( \frac{\log K}{\log x}\right) ^{2n}\right) =O\left( H\right) ,
\end{equation*}%
provided that, for example, 
\begin{equation*}
x=K^{\frac{\theta -1/2}{20n}}.
\end{equation*}

It remains to bound the fourth term. 
\begin{eqnarray*}
&&\frac{1}{H}\sum_{k=K}^{K+H-1}\left( \sigma _{x,g_{k}}-\frac{1}{2}\right)
^{2n}x^{2n\left( \sigma _{x,g_{k}}-\frac{1}{2}\right) }\left\{
\int_{1/2}^{\infty }x^{1/2-\sigma }\left| \sum_{p<x^{3}}\frac{\Lambda
_{x}\left( p\right) \log \left( xp\right) }{p^{\sigma +ig_{k}}}\right|
d\sigma \right\} ^{2n} \\
&\leq &\left\{ \frac{1}{H}\sum_{k=K}^{K+H-1}\left( \sigma _{x,g_{k}}-\frac{1%
}{2}\right) ^{4n}x^{4n\left( \sigma _{x,g_{k}}-\frac{1}{2}\right) }\right\}
^{1/2}\left\{ \frac{1}{H}\sum_{k=K}^{K+H-1}\left[ \int_{1/2}^{\infty
}x^{1/2-\sigma }\left| \sum_{p<x^{3}}\frac{\Lambda _{x}\left( p\right) \log
\left( xp\right) }{p^{\sigma +ig_{k}}}\right| d\sigma \right] ^{4n}\right\}
^{1/2}
\end{eqnarray*}%
by the Schwarz inequality. The first term in the product can be estimated as 
\begin{equation*}
O\left( \frac{1}{\left( \log x\right) ^{2n}}\right)
\end{equation*}%
by Lemma \ref{lemma_bound_integral} with $\xi =x^{4n}.$ For the second term,
we have 
\begin{eqnarray*}
&&\left[ \int_{1/2}^{\infty }x^{1/2-\sigma }\left| \sum_{p<x^{3}}\frac{%
\Lambda _{x}\left( p\right) \log \left( xp\right) }{p^{\sigma +ig_{k}}}%
\right| d\sigma \right] ^{4n} \\
&\leq &\left[ \int_{1/2}^{\infty }x^{1/2-\sigma }d\sigma \right]
^{4n-1}\int_{1/2}^{\infty }x^{1/2-\sigma }\left| \sum_{p<x^{3}}\frac{\Lambda
_{x}\left( p\right) \log \left( xp\right) }{p^{\sigma +ig_{k}}}\right|
^{4n}d\sigma \\
&=&\frac{1}{\left( \log x\right) ^{4n-1}}\int_{1/2}^{\infty }x^{1/2-\sigma
}\left| \sum_{p<x^{3}}\frac{\Lambda _{x}\left( p\right) \log \left(
xp\right) }{p^{\sigma +ig_{k}}}\right| ^{4n}d\sigma ,
\end{eqnarray*}%
where the second line follows by the H\"{o}lder inequality.

Hence, by Lemma \ref{lemma_bound_p_05}, we obtain 
\begin{eqnarray*}
&&\frac{1}{H}\sum_{k=K}^{K+H-1}\left( \sigma _{x,g_{k}}-\frac{1}{2}\right)
^{2n}x^{2n\left( \sigma _{x,g_{k}}-\frac{1}{2}\right) }\left\{
\int_{1/2}^{\infty }x^{1/2-\sigma }\left| \sum_{p<x^{3}}\frac{\Lambda
_{x}\left( p\right) \log \left( xp\right) }{p^{\sigma +ig_{k}}}\right|
d\sigma \right\} ^{2n} \\
&=&O\left( \sqrt{\log x}\int_{1/2}^{\infty }x^{1/2-\sigma }\frac{1}{H}%
\sum_{k=K}^{K+H-1}\left| \sum_{p<x^{3}}\frac{\Lambda _{x}\left( p\right)
\log \left( xp\right) }{p^{\sigma +ig_{k}}\log ^{2}x}\right| ^{4n}d\sigma
\right) \\
&=&O\left( \sqrt{\log x}\int_{1/2}^{\infty }x^{1/2-\sigma }d\sigma \right)
=O(1)
\end{eqnarray*}

provided that $x=K^{\frac{\theta -1/2}{20n}}.$ By using this in (\ref%
{Selberg_basic_estimate}), we find that 
\begin{equation*}
\sum_{k=K}^{K+H-1}\left| S\left( g_{k}\right) -S_{x}\left( g_{k}\right)
\right| ^{2n}=O\left( H\right) .
\end{equation*}%
$\square $

\section{Moments of the approximation to S(t)}

\label{section_moments_Sx}

In the next step we show that the moments of $S_{x}\left( t\right) $ are
approximately Gaussian.

\begin{lemma}
Suppose $1\leq c_{1}K^{\theta }\leq H\leq c_{2}K,$ where $\theta >1/2$ and $%
c_{1},c_{2}>0.$ Let $x\leq K^{\frac{2\theta -1}{6n}-\varepsilon }.$ Then,
for every integer $n\geq 1,$\newline
(i)%
\begin{equation*}
\sum_{k=K}^{K+H-1}\left| S_{x}\left( g_{k}\right) \right| ^{2n}=\frac{\left(
2n\right) !}{\left( 2\pi \right) ^{2n}n!}\left( H\left( \log \log K\right)
^{n}+O\left( H\left( \log \log K\right) ^{n-1}\right) \right) ,
\end{equation*}%
and (ii) 
\begin{equation*}
\sum_{k=K}^{K+H-1}S_{x}\left( g_{k}\right) ^{2n-1}=O(H).
\end{equation*}
\end{lemma}

\textbf{Proof: }First, we can write 
\begin{equation*}
S_{x}\left( t\right) =\frac{1}{2\pi i}\left( \eta -\overline{\eta }\right) ,
\end{equation*}%
where 
\begin{equation}
\eta =\eta \left( t\right) =\sum_{p<x^{3}}p^{-1/2-it}.
\label{definition_eta}
\end{equation}%
Hence, 
\begin{equation*}
\sum_{k=K}^{K+H-1}\left| S_{x}\left( g_{k}\right) \right| ^{2n}=\frac{1}{%
\left( 2\pi \right) ^{2n}}\sum_{l=0}^{2n}\left( -1\right) ^{n-l}\binom{2n}{l}%
\sum_{k=K}^{K+H-1}\eta \left( g_{k}\right) ^{l}\overline{\eta }\left(
g_{k}\right) ^{2n-l}.
\end{equation*}%
Here 
\begin{equation*}
\sum_{k=K}^{K+H-1}\eta \left( g_{k}\right) ^{l}\overline{\eta }\left(
g_{k}\right) ^{2n-l}=\sum_{p_{i}<x^{3}}\frac{1}{\sqrt{p_{1}\ldots p_{2n}}}%
\sum_{k=K}^{K+H-1}\left( \frac{p_{l+1}\ldots p_{2n}}{p_{1}\ldots p_{l}}%
\right) ^{ig_{k}}.
\end{equation*}

If $\left\{ p_{1},\ldots ,p_{l}\right\} \neq \left\{ p_{l+1},\ldots
,p_{2n}\right\} ,$ then by using Lemma \ref{lemma_estimate_exp_sum}, we
obtain 
\begin{equation*}
\sum_{k=K}^{K+H-1}\left( \frac{p_{l+1}\ldots p_{2n}}{p_{1}\ldots p_{l}}%
\right) ^{ig_{k}}\leq O\left( \left( Hn\frac{x^{n}\log x}{K^{1/2}\log K}%
+x^{n}K^{1/2}\log K\right) \right) .
\end{equation*}%
Since 
\begin{equation*}
\sum_{p_{i}<x^{3}}\frac{1}{\sqrt{p_{1}\ldots p_{2n}}}=O\left( \left(
\sum_{p\leq x^{3}}\frac{1}{\sqrt{p}}\right) ^{2n}\right) =O\left(
x^{2n}\right) ,
\end{equation*}%
hence in the case $l\neq n,$ we have 
\begin{eqnarray*}
\sum_{k=K}^{K+H-1}\eta \left( g_{k}\right) ^{l}\overline{\eta }\left(
g_{k}\right) ^{2n-l} &=&O\left( x^{2n}\left( Hn\frac{x^{n}\log x}{%
K^{1/2}\log K}+x^{n}K^{1/2}\log K\right) \right) \\
&=&O(H).
\end{eqnarray*}

If $l=n,$ then 
\begin{eqnarray*}
\sum_{k=K}^{K+H-1}\eta \left( g_{k}\right) ^{n}\overline{\eta }\left(
g_{k}\right) ^{n} &=&H\sum_{\substack{ p<x^{3}  \\ p_{1}\ldots
p_{n}=p_{n+1}\ldots p_{2n}}}\frac{1}{p_{1}\ldots p_{n}}+O(H) \\
&=&n!H\left( \sum_{p<x^{3}}\frac{1}{p}\right) ^{n}+O\left( n!H\sum_{p<x^{3}}%
\frac{1}{p_{1}\ldots p_{n-2}p_{n}^{2}}\right) +O\left( H\right) ,
\end{eqnarray*}%
where the second equality follows from the fact that the number of ways in
which a number of the form $p_{1}\ldots p_{n}$ can be written as a product
of $n$ primes is equal to $n!$ if the primes are all different and less than 
$n!$ if two or more of the primes are equal. Hence, 
\begin{eqnarray*}
\sum_{k=K}^{K+H-1}\eta \left( g_{k}\right) ^{n}\overline{\eta }\left(
g_{k}\right) ^{n} &=&n!H\left( \log \log x\right) ^{n}+O\left( n!H\left(
\log \log x\right) ^{n-1}\right) \\
&=&n!H\left( \log \log K\right) ^{n}+O\left( n!H\left( \log \log K\right)
^{n-1}\right) .
\end{eqnarray*}%
It follows that 
\begin{eqnarray*}
\sum_{k=K}^{K+H-1}\left| S_{x}\left( g_{k}\right) \right| ^{2n} &=&\frac{1}{%
\left( 2\pi \right) ^{2n}}\sum_{l=0}^{2n}\left( -1\right) ^{n-l}\binom{2n}{l}%
\sum_{k=K}^{K+H-1}\eta \left( g_{k}\right) ^{l}\overline{\eta }\left(
g_{k}\right) ^{2n-l} \\
&=&\frac{\left( 2n\right) !}{\left( 2\pi \right) ^{2n}n!}\left( H\left( \log
\log K\right) ^{n}+O\left( H\left( \log \log K\right) ^{n-1}\right) \right) .
\end{eqnarray*}%
The proof of (ii) is similar, except that in this case it is always true
that $\left\{ p_{1},\ldots ,p_{l}\right\} \neq \left\{ p_{l+1},\ldots
,p_{2n-1}\right\} .\square $

\section{Moments of S(t)}

\label{section_moments_S_and_proof_of_Thm1}

\begin{theorem}
Suppose $1\leq c_{1}K^{\theta }\leq H\leq c_{2}K,$ where $\theta >1/2$ and $%
c_{1},c_{2}>0.$ Then, for every integer $n\geq 1,$ \newline
(i)%
\begin{equation*}
\sum_{k=K}^{K+H-1}\left| S\left( g_{k}\right) \right| ^{2n}=\frac{\left(
2n\right) !}{\left( 2\pi \right) ^{2n}n!}\left( H\left( \log \log K\right)
^{n}+O\left( H\left( \log \log K\right) ^{n-1/2}\right) \right) ,
\end{equation*}%
and (ii) 
\begin{equation*}
\sum_{k=K}^{K+H-1}\left( S\left( g_{k}\right) \right) ^{2n-1}=O\left(
H\left( \log \log K\right) ^{n-1}\right) .
\end{equation*}
\end{theorem}

\textbf{Proof:} Take $x=K^{\frac{2\theta -1}{6n}-\varepsilon }.$ Then, the
triangle inequality for $L^{p}$ norms implies that%
\begin{eqnarray*}
\left| \left( \frac{1}{H}\sum_{k=K}^{K+H-1}\left| S\left( g_{k}\right)
\right| ^{2n}\right) ^{1/2n}-\left( \frac{1}{H}\sum_{k=K}^{K+H-1}\left|
S_{x}\left( g_{k}\right) \right| ^{2n}\right) ^{1/2n}\right| &\leq &\left( 
\frac{1}{H}\sum_{k=K}^{K+H-1}\left| S\left( g_{k}\right) -S_{x}\left(
g_{k}\right) \right| ^{2n}\right) ^{1/2n} \\
&=&O\left( 1\right) .
\end{eqnarray*}

Hence, 
\begin{equation*}
\left( \frac{1}{H}\sum_{k=K}^{K+H-1}\left| S\left( g_{k}\right) \right|
^{2n}\right) ^{1/2n}=\left( \frac{\left( 2n\right) !}{\left( 2\pi \right)
^{2n}n!}\left[ \left( \log \log K\right) ^{n}+O\left( \left( \log \log
K\right) ^{n-1}\right) \right] \right) ^{1/2n}+O\left( 1\right) ,
\end{equation*}%
and 
\begin{equation*}
\frac{1}{H}\sum_{k=K}^{K+H-1}\left| S\left( g_{k}\right) \right| ^{2n}=\frac{%
\left( 2n\right) !}{\left( 2\pi \right) ^{2n}n!}\left( \log \log K\right)
^{n}+O\left( \left( \log \log K\right) ^{n-1/2}\right) .
\end{equation*}%
For the proof of (ii), we estimate%
\begin{equation*}
S\left( g_{k}\right) ^{2n-1}-S_{x}\left( g_{k}\right) ^{2n-1}=O\left(
\sum_{\nu =1}^{2n-1}\left| S_{x}\left( g_{k}\right) \right| ^{2n-1-\nu
}\left| S\left( g_{k}\right) -S_{x}\left( g_{k}\right) \right| ^{\nu
}\right) ,
\end{equation*}%
and note that 
\begin{equation*}
\sum_{k=K}^{K+H-1}\left| S_{x}\left( g_{k}\right) \right| ^{2n-1-\nu }\left|
S\left( g_{k}\right) -S_{x}\left( g_{k}\right) \right| ^{\nu }\leq \left(
\sum_{k=K}^{K+H-1}\left| S_{x}\left( g_{k}\right) \right| ^{2n}\right) ^{%
\frac{2n-1-\nu }{2n}}\left( \sum_{k=K}^{K+H-1}\left| S\left( g_{k}\right)
-S_{x}\left( g_{k}\right) \right| ^{\frac{\nu }{1+\nu }2n}\right) ^{\frac{%
1+\nu }{2n}},
\end{equation*}%
where we used the H\"{o}lder inequality with $p=2n/\left( 2n-1-\nu \right) $
and $q=2n/\left( 1+\nu \right) .$ Next, we use the inequality 
\begin{equation*}
\left| \frac{1}{H}\sum_{k=K}^{K+H-1}\left| S\left( g_{k}\right) -S_{x}\left(
g_{k}\right) \right| ^{\frac{\nu }{1+\nu }2n}\right| ^{\frac{1+\nu }{\nu }%
\frac{1}{2n}}\leq \left| \frac{1}{H}\sum_{k=K}^{K+H-1}\left| S\left(
g_{k}\right) -S_{x}\left( g_{k}\right) \right| ^{2n}\right| ^{\frac{1}{2n}%
}=O\left( 1\right)
\end{equation*}%
in order to conclude that 
\begin{equation*}
\left| \sum_{k=K}^{K+H-1}\left| S\left( g_{k}\right) -S_{x}\left(
g_{k}\right) \right| ^{\frac{\nu }{1+\nu }2n}\right| ^{\frac{1+\nu }{2n}%
}\leq O\left( H^{\frac{1+\nu }{2n}}\right) .
\end{equation*}%
Therefore, 
\begin{eqnarray*}
\sum_{k=K}^{K+H-1}\left| S_{x}\left( g_{k}\right) \right| ^{2n-1-\nu }\left|
S\left( g_{k}\right) -S_{x}\left( g_{k}\right) \right| ^{\nu } &=&O\left( H^{%
\frac{2n-1-\nu }{2n}}\left( \log \log K\right) ^{n-\frac{1+\nu }{2}}H^{\frac{%
1+\nu }{2n}}\right) \\
&=&O\left( H\left( \log \log K\right) ^{n-\frac{1+\nu }{2}}\right) =O\left(
H\left( \log \log K\right) ^{n-1}\right)
\end{eqnarray*}%
for $1\leq \nu \leq 2n-1.$ Hence,%
\begin{eqnarray*}
\sum_{k=K}^{K+H-1}S\left( g_{k}\right) ^{2n-1}
&=&\sum_{k=K}^{K+H-1}S_{x}\left( g_{k}\right) ^{2n-1}+O\left( H\left( \log
\log K\right) ^{n-1}\right) \\
&=&O\left( H\left( \log \log K\right) ^{n-1}\right) .
\end{eqnarray*}

$\square $

\bigskip

\begin{corollary}
\label{corollary_moments}Suppose $1\leq c_{1}K^{\theta }\leq H\leq c_{2}K,$
where $\theta >1/2$ and $c_{1},c_{2}>0.$ Then, (i)
\end{corollary}

\begin{equation*}
\frac{1}{H}\sum_{k=K}^{K+H-1}\left| \frac{\sqrt{2}\pi S\left( t_{k}+\xi
\sigma _{k}\right) }{\sqrt{\log \log t_{k}}}\right| ^{2n}=\frac{\left(
2n\right) !}{\left( 2\pi \right) ^{2n}n!}\left( 1+O\left( \left( \log \log
K\right) ^{-1/2}\right) \right) ,
\end{equation*}%
and (ii) 
\begin{equation*}
\frac{1}{H}\sum_{k=K}^{K+H-1}\left( \frac{\sqrt{2}\pi S\left( t_{k}+\xi
\sigma _{k}\right) }{\sqrt{\log \log t_{k}}}\right) ^{2n-1}=O\left( \left(
\log \log K\right) ^{-1/2}\right) ,
\end{equation*}

This Corollary implies Theorem \ref{theorem_main_X} (and Theorem \ref%
{theorem_main} as a consequence).

\bigskip

\section{Covariance}

\label{section_covariance}

\begin{lemma}
\label{lemma_estimate_exp_sum_copy}Let two non-equal primes $p_{1},p_{2}$ be
both less than $y\leq cK$. Assume $1\leq H\leq cK,$ and let $k^{\prime
}=k+x, $ where $0<x<K^{\varepsilon },$ with $\varepsilon \in \lbrack 0,1).$
Then, 
\begin{equation*}
\sum_{k=K}^{K+H-1}\exp \left( -i\left( g_{k}\log p_{1}-g_{k^{\prime }}\log
p_{2}\right) \right) =O\left( H\frac{y^{1/2}\log y}{K^{1/2}\log K}%
+y^{1/2}K^{1/2}\log K\right) .
\end{equation*}
\end{lemma}

\textbf{Proof}: By using Lemma \ref{lemma_derivatives_g}, we can estimate: 
\begin{equation*}
g^{\prime \prime }\left( t\right) \log p_{1}-g^{\prime \prime }\left(
t+x\right) \log p_{2}=\frac{-2\pi }{t\left( \log t\right) ^{2}}\left( \log
p_{1}-\log p_{2}\right) +o\left( \frac{\log y}{t\left( \log t\right) ^{2}}%
\right) +o\left( \frac{x\log y}{t^{2}}\right) .
\end{equation*}%
It follows that 
\begin{equation*}
\left| g^{\prime \prime }\left( t\right) \log p_{1}-g^{\prime \prime }\left(
t+x\right) \log p_{2}\right| \geq c\frac{1}{yK\left( \log K\right) ^{2}},
\end{equation*}%
and 
\begin{equation*}
\left| g^{\prime \prime }\left( t\right) \log p_{1}-g^{\prime \prime }\left(
t+x\right) \log p_{2}\right| \leq cy\log y\frac{1}{yK\left( \log K\right)
^{2}}.
\end{equation*}%
The conclusion of the lemma follows by applying Theorem \ref%
{theorem_van_der_Corput}. $\square $

\begin{lemma}
\label{lemma_a_sum_over_primes} Suppose that $s\left( x\right) =c\left( \log
x\right) ^{\beta -1}+O\left( \left( \log x\right) ^{\beta -2}\right) ,$
where $\beta >0$. Then, 
\begin{equation*}
\sum_{p\leq x}\frac{1}{p}p^{is\left( x\right) }=\left( 1-\beta \right)
_{+}\log \log x+O\left( 1\right) .
\end{equation*}
\end{lemma}

\bigskip

\textbf{Proof:} This is a direct consequence of Lemma 3.4 in \cite%
{bourgade10}. $\square $

\textbf{Proof of Theorem \ref{theorem_covariance}}: In order to compute $%
\mathbb{E}X_{1}^{\left( N\right) }X_{2}^{\left( N\right) },$ we proceed as
above in the calculation of $\mathbb{E}\left( X_{1}^{\left( N\right)
}\right) ^{2}.$

Since by Proposition \ref{proposition_approximation}, 
\begin{equation*}
\sum_{k=K}^{K+H-1}\left| S\left( g_{k}\right) -S_{x}\left( g_{k}\right)
\right| ^{2n}=O\left( H\right) ,
\end{equation*}%
therefore, it is essential to compute 
\begin{equation*}
\frac{1}{N}\sum_{k=N}^{2N-1}S_{x}\left( g_{k}\right) S_{x}\left(
g_{k^{\prime }}\right) ,
\end{equation*}%
where $k^{\prime }=k+\left( \log N\right) ^{\beta }.$ By using function $%
\eta ,$ defined in (\ref{definition_eta}), we obtain: 
\begin{equation*}
\sum_{k=N}^{2N-1}S_{x}\left( g_{k}\right) S_{x}\left( g_{k^{\prime }}\right)
=-\frac{1}{\left( 2\pi \right) ^{2}}\sum_{k=N}^{2N-1}\left( \eta \left(
g_{k}\right) \eta \left( g_{k^{\prime }}\right) -\eta \left( g_{k}\right) 
\overline{\eta }\left( g_{k^{\prime }}\right) -\overline{\eta }\left(
g_{k}\right) \eta \left( g_{k^{\prime }}\right) +\overline{\eta }\left(
g_{k}\right) \overline{\eta }\left( g_{k^{\prime }}\right) \right) .
\end{equation*}%
For the first term in this sum, we write 
\begin{equation*}
\sum_{k=N}^{2N-1}\eta \left( g_{k}\right) \eta \left( g_{k^{\prime }}\right)
=\sum_{p_{1},p_{2}\leq x^{2}}\frac{1}{\sqrt{p_{1}p_{2}}}%
\sum_{k=N}^{2N-1}p_{1}^{-ig_{k}}p_{2}^{-ig_{k^{\prime }}}.
\end{equation*}%
Note that the sum 
\begin{equation*}
\sum_{k=N}^{2N-1}p_{1}^{-ig_{k}}p_{2}^{-ig_{k^{\prime
}}}=\sum_{k=N}^{2N-1}\exp [-i(g_{k}\log p_{1}+g_{k+\left( \log N\right)
^{\beta }}\log p_{2})]
\end{equation*}%
is an exponential sum, and it can be estimated by using van der Corput's
theorem by noticing that the second derivative of the function 
\begin{equation*}
g\left( s\right) \log p_{1}+g\left( s+\alpha \left( \log N\right) ^{\beta
}\right) \log p_{2}
\end{equation*}%
is bounded by $O\left( N^{-1}\left( \log N\right) ^{-2}\right) $ from below
and by $O\left( N^{-1}\left( \log N\right) ^{-2}\log x\right) $ from above.
This implies that with an appropriate choice of $x,$ 
\begin{equation*}
\frac{1}{N}\sum_{k=N}^{2N-1}\eta \left( g_{k}\right) \eta \left(
g_{k^{\prime }}\right) =O(1),
\end{equation*}%
and similarly for $N^{-1}\sum_{k=N}^{2N-1}\overline{\eta }\left(
g_{k}\right) \overline{\eta }\left( g_{k^{\prime }}\right) .$

Therefore 
\begin{equation*}
\frac{1}{N}\sum_{k=N}^{2N-1}S_{x}\left( g_{k}\right) S_{x}\left(
g_{k^{\prime }}\right) =\frac{1}{2\pi ^{2}}\mathrm{Re}\left[ \frac{1}{N}%
\sum_{k=N}^{2N-1}\eta \left( g_{k}\right) \overline{\eta }\left(
g_{k^{\prime }}\right) \right] +O\left( 1\right) ,
\end{equation*}%
where

\begin{equation*}
\sum_{k=N}^{2N-1}\eta \left( g_{k}\right) \overline{\eta }\left(
g_{k^{\prime }}\right) =\sum_{p_{1},p_{2}\leq x^{3}}\frac{1}{\sqrt{p_{1}p_{2}%
}}\sum_{k=N}^{2N-1}p_{1}^{-ig_{k}}p_{2}^{ig_{k^{\prime }}}.
\end{equation*}%
If $p_{1}\neq p_{2},$ then by using Lemma \ref{lemma_estimate_exp_sum_copy},
we can estimate 
\begin{equation*}
\sum_{p_{1},p_{2}\leq x^{2}}\frac{1}{\sqrt{p_{1}p_{2}}}%
\sum_{k=N}^{2N-1}p_{1}^{-ig_{k}}p_{2}^{ig_{k^{\prime }}}=O\left(
x^{3}N^{1/2}\log N\right) =O(N),
\end{equation*}%
provided that $x=N^{\kappa }$ and $\kappa \leq 1/6.$

If $p_{1}=p_{2},$ then we have 
\begin{equation*}
\frac{1}{N}\sum_{k=N}^{2N-1}\sum_{p\leq x^{3}}\frac{1}{p}p^{-i\left(
g_{k}-g_{k^{\prime }}\right) }
\end{equation*}

If one sets $x=N^{\kappa },$ then by using the definition of function $g,$
it is easy to see that for every $k\in \left[ N,2N-1\right] ,$ and $%
k^{\prime }=k+\left( \log N\right) ^{\beta },$ we have%
\begin{eqnarray*}
g_{k^{\prime }}-g_{k} &=&2\pi (\log N)^{\beta -1}+O\left( (\log N)^{\beta
-2}\right) \\
&=&\left( 2\pi /\kappa ^{\beta -1}\right) (\log x)^{\beta -1}+O\left( (\log
x)^{\beta -2}\right) ,
\end{eqnarray*}%
where the implicit constant in the $O$-term does not depend on $k.$

Hence, by Lemma \ref{lemma_a_sum_over_primes}, 
\begin{equation*}
\frac{1}{N}\sum_{k=N}^{2N-1}\sum_{p\leq x^{3}}\frac{1}{p}p^{-i\left(
g_{k}-g_{k^{\prime }}\right) }=\left( 1-\beta \right) _{+}\log \log
x+O\left( 1\right) .
\end{equation*}

It follows that 
\begin{equation*}
\frac{1}{N}\sum_{k=N}^{2N-1}S_{x}\left( g_{k}\right) S_{x}\left(
g_{k^{\prime }}\right) =\frac{1}{2\pi ^{2}}\left( 1-\beta \right) _{+}\log
\log x+O\left( 1\right) .
\end{equation*}

Next we note that 
\begin{eqnarray*}
\frac{1}{N}\sum_{k=N}^{2N-1}S\left( g_{k}\right) S\left( g_{k^{\prime
}}\right) &=&\frac{1}{N}\sum_{k=N}^{2N-1}S_{x}\left( g_{k}\right)
S_{x}\left( g_{k^{\prime }}\right) \\
&&+\frac{1}{N}\sum_{k=N}^{2N-1}S_{x}\left( g_{k}\right) \left( S\left(
g_{k^{\prime }}\right) -S_{x}\left( g_{k^{\prime }}\right) \right) \\
&&+\frac{1}{N}\sum_{k=N}^{2N-1}\left( S\left( g_{k}\right) -S_{x}\left(
g_{k}\right) \right) S_{x}\left( g_{k^{\prime }}\right) \\
&&+\frac{1}{N}\sum_{k=N}^{2N-1}\left( S\left( g_{k}\right) -S_{x}\left(
g_{k}\right) \right) \left( S\left( g_{k^{\prime }}\right) -S_{x}\left(
g_{k^{\prime }}\right) \right) .
\end{eqnarray*}%
By the Schwarz inequality, the last three terms can be estimated as $\left(
\log \log N\right) ^{1/2},$ and therefore we have 
\begin{equation*}
\frac{1}{N}\sum_{k=N}^{2N-1}S\left( g_{k}\right) S\left( g_{k^{\prime
}}\right) =\frac{1}{2\pi ^{2}}\left( 1-\beta \right) _{+}\log \log x+O\left(
\left( \log \log x\right) ^{1/2}\right) ,
\end{equation*}

This implies that 
\begin{equation*}
\mathbb{E}X_{1}^{\left( N\right) }X_{2}^{\left( N\right) }=\left( 1-\beta
\right) _{+}.
\end{equation*}%
$\square $

\section{Joint Moments}

\label{section_joint_moments}

\textbf{Proof of Theorem \ref{theorem_joint_moments_X}:} It is clearly
enough to prove the corresponding result for random variables $S\left(
g_{k_{1}}\right) $ and $S(g_{k_{2}})$ since $X_{i}^{\left( N\right) }$ are
the rescaled versions of these random variables. In fact, as a consequence
of the Selberg approximation result, it is enough to show that $S_{x}\left(
g_{k_{1}}\right) \,$and $S_{x}\left( g_{k_{2}}\right) $ have the required
moments.

Indeed, 
\begin{eqnarray*}
S\left( g_{k_{1}}\right) ^{a}S\left( g_{k_{2}}\right) ^{b} &=&\left(
S_{x}\left( g_{k_{1}}\right) +S\left( g_{k_{1}}\right) -S_{x}\left(
g_{k_{1}}\right) \right) ^{a}\left( S_{x}\left( g_{k_{2}}\right) +S\left(
g_{k_{2}}\right) -S_{x}\left( g_{k_{2}}\right) \right) ^{b} \\
&=&S_{x}\left( g_{k_{1}}\right) ^{a}S_{x}\left( g_{k_{2}}\right) ^{b} \\
&&+O\left( \sum_{s,t}S_{x}\left( g_{k_{1}}\right) ^{s}\left( S\left(
g_{k_{1}}\right) -S_{x}\left( g_{k_{1}}\right) \right) ^{a-s}S_{x}\left(
g_{k_{2}}\right) ^{t}\left( S\left( g_{k_{2}}\right) -S_{x}\left(
g_{k_{2}}\right) \right) ^{b-t}\right) ,
\end{eqnarray*}%
where the sum is over $s$ and $t$ such that $0\leq s\leq a,$ $0\leq t\leq b,$
and $s+t < a+b.$

After we sum over $k_{1}$ and apply the Schwarz inequality twice, we find
that 
\begin{eqnarray*}
&&\frac{1}{N}\sum_{k_{1}=N}^{2N-1}\left( S\left( g_{k_{1}}\right)
^{a}S\left( g_{k_{2}}\right) ^{b}-S_{x}\left( g_{k_{1}}\right)
^{a}S_{x}\left( g_{k_{2}}\right) ^{b}\right) \\
&=&O\left( \sum_{s,t}\left( \frac{1}{N}\sum_{k_{1}=N}^{2N-1}S_{x}\left(
g_{k_{1}}\right) ^{4s}\right) ^{1/4}\left( \frac{1}{N}%
\sum_{k_{1}=N}^{2N-1}S_{x}\left( g_{k_{2}}\right) ^{4t}\right) ^{1/4}\right.
\\
&&\left. \times \left( \frac{1}{N}\sum_{k_{1}=N}^{2N-1}\left( S\left(
g_{k_{1}}\right) -S_{x}\left( g_{k_{1}}\right) \right) ^{4(a-s)}\right)
^{1/4}\left( \frac{1}{N}\sum_{k_{1}=N}^{2N-1}\left( S\left( g_{k_{1}}\right)
-S_{x}\left( g_{k_{1}}\right) \right) ^{4(a-s)}\right) ^{1/4}\right) \\
&=&O\left( \left( \log \log N\right) ^{(a+b-1)/2}\right) .
\end{eqnarray*}%
Hence, if variables $S$ and $S_{x}$ are scaled by $\left( \log \log N\right)
^{-1},$ the difference in their moments is of order $\left( \log \log
N\right) ^{-1/2}.$

The result about moments of the scaled versions of $S_{x}\left(
g_{k_{1}}\right) \,$and $S_{x}\left( g_{k_{2}}\right) $ follows from the
result for random variables 
\begin{equation*}
\eta _{i}^{\left( N\right) }:=\frac{1}{\sqrt{\log \log N}}\eta \left(
g_{k_{i}\left( N,\omega \right) }\right) ,
\end{equation*}%
where $i=1,2,$ and $\eta \left( t\right) $ is as defined in (\ref%
{definition_eta}).

\begin{theorem}
\label{theorem_joint_moments_eta}Let $a_{1},a_{2},b_{1},b_{2}\geq 0.$ The
joint moments \ of random variables $\eta _{1}^{\left( N\right) }$ and $\eta
_{2}^{\left( N\right) }$,%
\begin{equation*}
m_{N}\left( a_{1},a_{2},b_{1},b_{2}\right) :=\mathbb{E}\left( \eta
_{1}^{\left( N\right) }\right) ^{a_{1}}\left( \overline{\eta _{1}^{\left(
N\right) }}\right) ^{a_{2}}\left( \eta _{2}^{\left( N\right) }\right)
^{b_{1}}\left( \overline{\eta _{2}^{\left( N\right) }}\right) ^{b_{2}},
\end{equation*}%
converge to the corresponding joint moments of complex Gaussian random
variables $\eta _{1}$ and $\eta _{2},$ which have the following covariance
structure: $\mathbb{E\eta }_{i}^{2}=\mathbb{E}\overline{\mathbb{\eta }}%
_{i}^{2}=\mathbb{E\eta }_{1}\eta _{2}=\mathbb{E}\overline{\mathbb{\eta }}_{1}%
\overline{\eta }_{2}=0,$ $\mathbb{E\eta }_{i}\overline{\eta _{i}}=1,$ $%
\mathbb{E\eta }_{1}\overline{\eta }_{2}=\mathbb{E}\overline{\mathbb{\eta }}%
_{1}\eta _{2}=\left( 1-\beta \right) _{+}.$
\end{theorem}

Indeed, if this result holds, then the joint moments of (real) random
variables 
\begin{equation*}
S_{x}\left( g_{k_{j}}\right) =\frac{1}{2\pi i}\left( \eta _{j}^{\left(
N\right) }-\overline{\eta }_{j}^{(N)}\right)
\end{equation*}%
converge to the corresponding joint moments of Gaussian random variables $%
S_{1}$ and $S_{2},$ where 
\begin{equation*}
E\left( S_{1}^{2}\right) =E\left( S_{2}^{2}\right) =\frac{1}{2\pi ^{2}},
\end{equation*}%
and 
\begin{equation*}
E\left( S_{1}S_{2}\right) =\frac{1}{2\pi ^{2}}\left( 1-\beta \right) _{+}.
\end{equation*}%
This implies the statement of Theorem \ref{theorem_joint_moments_X}.

Before attacking Theorem \ref{theorem_joint_moments_eta}, let us recall the
Wick Rule for the joint moments of Gaussian random variables, namely, 
\begin{equation*}
\mathbb{E}\left[ x_{i_{1}}\ldots x_{i_{k}}\right] =\sum_{\pi \in P_{2}\left(
\left\{ 1,\ldots ,k\right\} \right) }\prod_{\left( r,s\right) \in \pi }%
\mathbb{E}\left[ x_{i_{r}}x_{i_{s}}\right] ,
\end{equation*}%
where the sum is over all pairings of indices $1,\ldots ,k.$ (In particular,
if $k$ is odd, then the sum is empty.) (See, for example, Theorem 22.3 in %
\cite{nica_speicher06} or Appendix 1 on p. 13 in \cite{zee03}).

If we apply this rule to random variables $\eta _{i},\overline{\eta }_{i},$
then we find that 
\begin{equation*}
m\left( a_{1},a_{2},b_{1},b_{2}\right) :=\mathbb{E}\left( \eta _{1}\right)
^{a_{1}}\left( \overline{\eta _{1}}\right) ^{a_{2}}\left( \eta _{2}\right)
^{b_{1}}\left( \overline{\eta _{2}}\right) ^{b_{2}}
\end{equation*}%
is zero unless $a_{1}+b_{1}=a_{2}+b_{2}.$ If $a_{1}+b_{1}=a_{2}+b_{2},$ then 
\begin{equation}
m\left( a_{1},a_{2},b_{1},b_{2}\right) =n\left(
k,a_{1},a_{2},b_{1},b_{2}\right) \left( 1-\beta \right) _{+}^{k},
\label{formula_joint_moment_Gaussian}
\end{equation}%
where $n\left( k,a_{1},a_{2},b_{1},b_{2}\right) $ is a number of ways to
pair $a_{1}$ elements $\eta _{1}$ and $b_{1}$ elements $\eta _{2}$ with $%
a_{2}$ elements $\overline{\eta }_{1}$ and $b_{2}$ elements $\overline{\eta }%
_{2}$ so that exactly $k$ elements are connected with an element that has a
different index.

Also, we need a generalization of Lemma \ref{lemma_estimate_exp_sum}.

\begin{lemma}
\label{lemma_estimate_exp_sum_copy2}Let 
\begin{eqnarray*}
h\left( x\right) &=&g\left( x\right) \left( \sum_{k=1}^{a_{2}}\log
q_{k}-\sum_{k=1}^{a_{1}}\log p_{k}\right) \\
&&+g\left( x+u\right) \left( \sum_{k=a_{2}+1}^{a_{2}+b_{2}}\log
q_{k}-\sum_{k=a_{1}+1}^{a_{1}+b_{1}}\log p_{k}\right) ,
\end{eqnarray*}%
where $g\left( x\right) $ is as defined in (\ref{definition_g}), $\left\{
p_{1},\ldots ,p_{a_{1}+b_{1}}\right\} \neq \left\{ q_{1},\ldots
q_{a_{2}+b_{2}}\right\} $ and primes $p_{i}$ and $q_{i}$ are less than $%
y\leq K^{\varepsilon }$ for all $i$ and $\varepsilon <1/n$. Assume $1\leq
H\leq cK$ and $u\leq \alpha (\log K)^{\beta }.$ Let $n=\left\lfloor \left(
a_{1}+a_{2}+b_{1}+b_{2}\right) /2\right\rfloor .$ Then, 
\begin{equation*}
\sum_{k=K}^{K+H-1}e^{ih\left( k\right) }=O\left( H\frac{y^{n/2}\log y}{%
K^{1/2}\log K}+y^{n/2}K^{1/2}\log K\right) .
\end{equation*}
\end{lemma}

\textbf{Proof:} We can re-write the definition of $h\left( x\right) $ as
follows: 
\begin{eqnarray*}
h\left( x\right) &=&h_{1}\left( x\right) +h_{2}\left( x\right) =g\left(
x\right) \left( \sum_{k=1}^{a_{2}+b_{2}}\log
q_{k}-\sum_{k=1}^{a_{1}+b_{1}}\log p_{k}\right) \\
&&+\left( g\left( x+u\right) -g\left( x\right) \right) \left(
\sum_{k=a_{2}+1}^{a_{2}+b_{2}}\log q_{k}-\sum_{k=a_{1}+1}^{a_{1}+b_{1}}\log
p_{k}\right) .
\end{eqnarray*}%
The second derivative of the first term can be estimated as in Lemma \ref%
{lemma_estimate_exp_sum}: If $x\in \left[ K,K+H\right] ,$ then 
\begin{equation*}
h_{1}^{\prime \prime }\left( x\right) \in \left[ \lambda ,\kappa \lambda %
\right] ,
\end{equation*}%
where 
\begin{equation*}
\lambda =\frac{c_{1}}{y^{n}K\left( \log K\right) ^{2}}\text{ and }\kappa
=c_{2}y^{n}\log y.
\end{equation*}

For the second term, we note that 
\begin{equation*}
\left( g\left( x+u\right) -g\left( x\right) \right) ^{\prime \prime
}=g^{\prime \prime \prime }\left( \theta \right) u,
\end{equation*}%
where $\theta \in \left[ x,x+u\right] ,$ and \ by using Lemma \ref%
{lemma_derivatives_g} we find that 
\begin{equation*}
h_{2}^{\prime \prime }\left( x\right) =O\left( \frac{\left( \log K\right)
^{\beta -2}}{K^{2}}\log y\right) =o(h_{1}^{\prime \prime }\left( x\right) ,
\end{equation*}%
provided that $y\leq K^{\varepsilon }$ with $\varepsilon <1/n.$ It follows
that $h^{\prime \prime }\left( x\right) \sim h_{1}^{\prime \prime }\left(
x\right) $, and the conclusion of the lemma follows by an application of
Theorem \ref{theorem_van_der_Corput} as in Lemma \ref{lemma_estimate_exp_sum}%
. $\square $

\textbf{Proof of Theorem \ref{theorem_joint_moments_eta}:} By definition, we
write 
\begin{eqnarray*}
\mathbb{E}\left( \eta _{1}^{\left( N\right) }\right) ^{a_{1}}\left( 
\overline{\eta _{1}^{\left( N\right) }}\right) ^{a_{2}}\left( \eta
_{2}^{\left( N\right) }\right) ^{b_{1}}\left( \overline{\eta _{2}^{\left(
N\right) }}\right) ^{b_{2}} &=&\frac{1}{N\left( \log \log N\right) ^{\left(
a_{1}+a_{2}+b_{1}+b_{2}\right) /2}} \\
&&\times \sum_{k_{1}=N}^{2N-1}\left( \sum_{p\leq x^{2}}\frac{p^{-ig_{k_{1}}}%
}{\sqrt{p}}\right) ^{a_{1}}\left( \sum_{q\leq x^{2}}\frac{q^{ig_{k_{1}}}}{%
\sqrt{q}}\right) ^{a_{2}} \\
&&\times \left( \sum_{p\leq x^{2}}\frac{p^{-ig_{k_{2}}}}{\sqrt{p}}\right)
^{b_{1}}\left( \sum_{q\leq x^{2}}\frac{q^{ig_{k_{2}}}}{\sqrt{q}}\right)
^{b_{2}},
\end{eqnarray*}%
where $k_{2}=k_{1}+\left[ \alpha \left( \log N\right) ^{\beta }\right] .$ If
we expand the product of sums, we get for a general term 
\begin{equation*}
t(p,p^{\prime },q,q^{\prime }):=\frac{1}{\sqrt{p_{1}\ldots
p_{a_{1}+b_{1}}q_{1}\ldots q_{a_{2+b_{2}}}}}\frac{\left( q_{1}\ldots
q_{a_{2}}\right) ^{ig_{k_{1}}}\left( q_{a_{2}+1}\ldots
q_{a_{2}+b_{2}}\right) ^{ig_{k_{2}}}}{\left( p_{1}\ldots p_{a_{1}}\right)
^{ig_{k_{1}}}\left( p_{a_{1}+1}\ldots p_{a_{1}+b_{1}}\right) ^{ig_{k_{2}}}},
\end{equation*}%
where $p:=\left( p_{1},\ldots ,p_{a_{1}}\right) ,$ $p^{\prime
}:=(p_{a_{1}+1,}\ldots ,p_{a_{1}+a_{2}}),$ $q:=\left( q_{1},\ldots
,q_{b_{1}}\right) ,$ and $q^{\prime }:=(q_{b_{1}+1,}\ldots
,q_{b_{1}+b_{2}}). $

By using Lemma \ref{lemma_estimate_exp_sum_copy2}, we find that after we sum
this term over $k_{1}$ and divide it by \newline
$N\left( \log \log N\right) ^{\left( a_{1}+a_{2}+b_{1}+b_{2}\right) /2},$ we
get a non-negligible contribution if and only if there is a pairing that
puts every $q_{i}$ in a correspondence with a $p_{j},$ so that $q_{i}=p_{j}.$
In particular, it must be true that $a_{1}+b_{1}=a_{2}+b_{2}=n.$ Hence, the
moment is asymptotically equivalent to 
\begin{equation}
\frac{1}{N\left( \log \log N\right) ^{n}}\sum_{k_{1}=N}^{2N-1}\sum_{p\cdot
p^{\prime }=q\cdot q^{\prime }}t(p,p^{\prime },q,q^{\prime }),
\label{sum_joint_moment}
\end{equation}%
where $p\cdot p^{\prime }$ denotes the product of primes in $p$ and $%
p^{\prime }$, and similar for $q\cdot q^{\prime }.$ If we consider the sum
over all $\left( p,p^{\prime }\right) ,$ in which at least one $p_{i}$
appears twice, then we can see that this sum can be estimated as 
\begin{equation*}
O\left( \sum_{k_{1}=N}^{2N-1}\left( \sum_{p\leq x^{2}}\frac{1}{p^{2}}\right)
\left( \sum_{p\leq x^{2}}\frac{1}{p}\right) ^{n-2}\right) =O(N\left( \log
\log N\right) ^{n-2}),
\end{equation*}%
which gives a negligible contribution to the moment.

Otherwise, if every prime appears only once in $\left( p,p^{\prime }\right) $
and if $p\cdot p^{\prime }=q\cdot q^{\prime },$ then there is a unique
pairing between elements of $\left( p,p^{\prime }\right) $ and $\left(
q,q^{\prime }\right) .$ Let this pairing be called $\pi .$ That is, $\pi
\left( i\right) =j$ means that $p_{i}=q_{j}.$

The terms that satisfy pairing $\pi $ give the following contribution to the
sum in (\ref{sum_joint_moment}): 
\begin{equation*}
\sum_{k_{1}=N}^{2N-1}\left( \sum_{p\leq x^{2}}\frac{1}{p}\right)
^{n_{11}+n_{22}}\left( \sum_{p\leq x^{2}}\frac{p^{i\left(
g_{k_{1}}-g_{k_{2}}\right) }}{p}\right) ^{n_{12}}\left( \sum_{p\leq x^{2}}%
\frac{p^{-i\left( g_{k_{1}}-g_{k_{2}}\right) }}{p}\right)
^{n_{21}}+O(N\left( \log \log N\right) ^{n-2}),
\end{equation*}%
where 
\begin{eqnarray*}
n_{11} &=&\left\vert \left\{ i,j:\pi \left( i\right) =j,1\leq i\leq
a_{2},1\leq j\leq a_{1}\right\} \right\vert , \\
n_{12} &=&\left\vert \left\{ i,j:\pi \left( i\right) =j,1\leq i\leq
a_{2},a_{1}\leq j\leq a_{1}+b_{1}\right\} \right\vert ,
\end{eqnarray*}%
and so on.

By Lemma \ref{lemma_a_sum_over_primes}, this can be computed as 
\begin{equation*}
\left( \left( 1-\beta \right) _{+}\right) ^{n_{12}+n_{21}}N\left( \log \log
N\right) ^{n}+O(N\left( \log \log N\right) ^{n-1}).
\end{equation*}

After summing over all pairings we find that the moment equals 
\begin{equation*}
\sum_{\pi }\left( \left( 1-\beta \right) _{+}\right)
^{n_{12}+n_{21}}+O\left( \left( \log \log N\right) ^{-1}\right) .
\end{equation*}

Recall $n\left( k,a_{1},a_{2},b_{1},b_{2}\right) $ is a number of ways to
pair $a_{1}$ elements $\eta _{1}$ and $b_{1}$ elements $\eta _{2}$ with $%
a_{2}$ elements $\overline{\eta }_{1}$ and $b_{2}$ elements $\overline{\eta }%
_{2}$ so that exactly $k$ elements are connected with an element that has a
different index. That is, $n\left( k,a_{1},a_{2},b_{1},b_{2}\right) $ is the
number of pairings $\pi $ for which $n_{12}+n_{21}=k.$

It follows that asymptotically, the moment tends to $n\left(
k,a_{1},a_{2},b_{1},b_{2}\right) \left( 1-\beta \right) _{+}^{k},$ which is
exactly the corresponding joint moment of the Gaussian variables that we
obtained in formula (\ref{formula_joint_moment_Gaussian}). $\square $

\section{Conclusion}

We have shown that the distribution of two zeta zero fluctuations $f_{k}$
and $f_{k+x}$ approaches a two-variate Gaussian distribution with covariance 
$\left( 1-\beta \right) _{+}$, provided that $x\sim \left( \log k\right)
^{\beta }$. This gives an analogue of Gustavsson's results for fluctuations
of eigenvalues from Gaussian Unitary Ensemble. It is of obvious interest to
study the correlation of zeros at shorter distances. However, methods of
this paper are not easy to generalize to this case.

Some of the methods in this paper could perhaps be useful to extend
Gustavsson's results to other ensembles of random matrices, in particular to
the ensemble of uniformly distributed unitary random matrices. The proof
would proceed along the similar lines by using the additional tool by
Diaconis and Shashahani (\cite{diaconis_shahshahani94}) about expected
values of traces of moments of $U.$

\bigskip

\label{section_conclusion}

\bibliographystyle{plain}
\bibliography{comtest}
\qquad

\end{document}